\UseRawInputEncoding
\documentclass{amsart}
\usepackage{amsmath}
\usepackage{amssymb}
\newtheorem{thm}{Theorem}[section]
\newtheorem{cor}[thm]{Corollary}
\newtheorem{lem}[thm]{Lemma}
\newtheorem{prop}[thm]{Proposition}
\newtheorem{cons}[thm]{Construction}
\newtheorem{ex}[thm]{Example}
\theoremstyle{definition}
\newtheorem{defn}[thm]{Definition}
\newtheorem{notn}[thm]{Notation}
\theoremstyle{remark}
\newtheorem{rem}[thm]{Remark}

\long\def\Cor#1{\begin{cor} #1 \end{cor}}
\long\def\Lem#1{\begin{lem} #1 \end{lem}}
\long\def\Prop#1{\begin{prop} #1 \end{prop}}

\long\def\Def#1{\begin{defn} #1 \end{defn}}

\long\def\Rem#1{\begin{rem} #1 \end{rem}}
\long\def\Ex#1{\begin{ex} #1 \end{ex}}

\def\Sect{\section}
\def\Rarr#1#2{\xrightarrow[#2]{#1}}

\def\Darr#1#2{{\scriptstyle #1}\downarrow{\scriptstyle #2}}
\def\Uarr#1#2{{\scriptstyle #1}\uparrow{\scriptstyle #2}}
\long\def\Ref#1#2#3#4#5#6{
\bibitem{#1}
{\rm #2,}
\textit{#3.}
{\rm #4}
\textbf{#5}
{\rm #6.}
}
\long\def\Refb#1#2#3#4{
\bibitem{#1}
{\rm #2,}
\textit{#3.}
#4.
}
\def\Zz{{\mathbb Z}}
\def\Rr{{\mathbb R}}
\def\Cc{{\mathbb C}}
\def\Tt{{\mathbb T}}
\def\Ff{{\mathbb F}}

\def\into{\hookrightarrow}
\def\iso{\cong}
\def\leq{\leqslant}
\def\geq{\geqslant}

\def\comp{\mathbin{\mathchoice
{\circ}
{{\scriptstyle\circ}}
{{\scriptscriptstyle\circ}}
{{\scriptscriptstyle\circ}}
}}
\def\st{\mid}

\def\phi{\varphi}

\def\supp{{\rm supp}}

\def\KO{K{\rm O}}
\def\ind{{\rm ind}}
\def\Sq{{\rm Sq}}
\begin{document}

\title{New perspectives on a classical embedding theorem}
\author{M.~C.~Crabb}
\address{%
Institute of Mathematics\\
University of Aberdeen \\
Aberdeen AB24 3UE, UK}
\email{m.crabb@abdn.ac.uk}
\date{April 2025}
\begin{abstract}
In this expository note, recent results of Kishimoto and  Matsushita 
on triangulated manifolds
are linked to the classical criterion on the normal Stiefel-Whitney 
classes for existence of
an embedding of a smooth closed manifold into Euclidean space of
given dimension.
We also look back at Atiyah's K-theoretic condition for the existence
of a smooth embedding.
\end{abstract}

\maketitle
\Sect{Introduction}
Consider a closed smooth manifold of dimension $d$, with tangent bundle
$\tau M$. If $M$ embeds as a submanifold in a Euclidean space $\Rr^m$
with normal bundle $\nu$ of dimension $m-d$, it is a result from the
early days of Algebraic Topology,
attributed in \cite{massey} to Seifert and Whitney,
that the Stiefel-Whitney classes $w_i(\nu )=w_i(-\tau M)$ are zero
for $i\geq m-d$.

Another early result from the 1930s due to van Kampen and Flores 
\cite{kampen, flores} asserts that, for the standard
triangulation $K$ of the sphere $S^d$ as the boundary of
a $(d+1)$-simplex, if $d=2q+1$ is odd and $f : K^{(q)}\to
\Rr^{2q+1}$ is a continuous map on the $q$-skeleton,
then there exist two points $x,\, y\in K^{(q)}$, lying
in disjoint simplices, such that $f(x)=f(y)$.

These two classical results, on the topological embedding of a 
smooth manifold in Euclidean space and the combinatorial 
van Kampen-Flores theorem, are combined 
in a recent paper \cite{kish} by Kishimoto and Matsushita.

\medskip

This note is purely expository and is written primarily
as an account of the work of Kishimoto and Matsushita. We
shall also make connections to a paper of Frick and Harrison
\cite{frick3}.
And an appendix looks back at Atiyah's early
application of topological $K$-theory \cite{mfa}
to embeddings (and immersions).

\medskip

Throughout the paper, $M$ will be a closed smooth manifold
of dimension $d$, and in Section 3 the manifold 
$M$ will be triangulated
as a finite simplicial complex $K$.
The following notation will also be used without specific comment.
The cyclic group of order $2$ is written as $G=\Zz /2$, 
We write $L$ for the representation $\Rr$ of
$G$ with the action of the generator as $-1$.
For a free compact $G$-ENR $X$, we shall employ the systematic notation
$Y=X/G$ for the orbit
space and $\lambda$ for the line bundle $X\times_G L$ over $Y$.
We use cohomology with $\Ff_2$-coefficients throughout this
note (except that in the Appendix we also need $\Zz$-cohomology), writing simply $H^*$ for cohomology groups
and $H^*_G$ for the
equivariant Borel cohomology with $\Ff_2$-coefficients.
The $\Ff_2$-Euler class of $\lambda$,
$e(\lambda )\in H^1(Y)$, corresponds to the equivariant
cohomology Euler class $e(L)\in H^1_G(X)$ of the trivial bundle
$X\times L$.
\Sect{Smooth manifolds: the classical embedding theorem}
We begin with a result on smooth manifolds.
\Lem{\label{embed}
Suppose that $M$ is a closed smooth manifold of dimension $d$
with tangent bundle $\tau M$.
Consider the product $M\times M$ with the action of $G=\Zz /2$ that
interchanges the factors. 
We write $D(L\otimes \tau M)$ for the closed unit disc bundle
in $L\otimes\tau M$ for a chosen inner product on $\tau M$,
and $B(L\otimes\tau M)$ for the open unit disc bundle.
Choose an equivariant tubular neighbourhood
$\iota :D(L\otimes\tau M)\into M\times M$ of the diagonal and
take $X=M\times M-B(L\otimes\tau M)$.
Then, for an integer $m\geq 1$,
$e(\lambda )^m$ is non-zero if $w_i(-\tau M)\not=0$
for some $i\geq m-d$.
}
\begin{proof}
Suppose that $e(\lambda )^m=0$. We shall show that $w_i(-\tau M)=0$ for $i\geq m-d$.

Writing $\Delta : M\to M\times M$ for the diagonal inclusion
(and $\Delta_!$ for the associated Umkehr or direct image 
homomorphism),
we have a commutative diagram in $\Ff_2$-cohomology:
$$
\begin{matrix}
H^{m-d}(M) &\Rarr{\Delta_!}{\into}& H^m(M\times M) 
&\to& H^m(M\times M-B(\tau M))\\
\Uarr{}{}&&\Uarr{}{}&& \Uarr{}{}\\
H^{m-d}_G(M) &\Rarr{\Delta_!}{} &H^m_G (M\times M)&\to&
H^m_G(M\times M-B(L\otimes\tau M))\\
\Darr{=}{\phantom{=}}&&\Darr{\Delta^*}{\phantom{\Delta^*}}&&\Darr{\iota^*}{\phantom{\iota^*}}\\
H^{m-d}_G(M)&\Rarr{\cdot \,e(L\otimes\tau M)}{\into}&
H^{m}_G(M)&\to& H^m_G(S(L\otimes\tau M))
\end{matrix}
$$
In the diagram, the first row is the cohomology exact sequence of the
pair $(M\times M,$ $M\times M-B(\tau M))$;
the relative group $H^m(M\times M,M\times M-B(\tau M))$ is
identified by excision with $H^m(D(\tau M),S(\tau M))$
and then, via the Thom isomorphism, with $H^{m-d}(M)$.
The second row is similarly the $G$-equivariant Borel cohomology
sequence for the pair $(M\times M, M\times M-B(L\otimes\tau M))$
and maps to the first row by forgetting the group action.
The third row is the Gysin exact sequence of the sphere bundle
$S(L\otimes \tau M)$ over $M$, realized as the  
exact sequence
of the pair $(D(L\otimes\tau M), S(L\otimes\tau M))$.
The maps from the second to the third row are induced by the inclusion
$\iota : D(L\otimes\tau M)\into M\times M$.

Since $e(\lambda )^m=0$, there is a class 
$x\in H^{m-d}_G(M)$ such that $\Delta_!(x)=e(L)^m$
in the second row, and $x\cdot e(L\otimes \tau M)=e(L)^m$.
The non-equivariant top line,
$\Delta_! :H^{m-d}(M)\to H^m(M\times M)$
is injective, split by $\pi_!$, where $\pi :M\times M\to M$
is the projection to the first factor, and so
$x$ restricts to zero in the Gysin sequence of $S(L)$:
$$
\cdots \to H^{m-d-1}_G(M) \Rarr{\cdot e(L)}{}
H^{m-d}_G(M)\to H^{m-d}(M)\to \cdots \,.
$$
Thus $x=e(L)\cdot a$ for some $a\in H^{m-d-1}_G(M)$
and we obtain the identity
$$
a\cdot e(L)\cdot e(L\otimes\tau M) =e(L)^m
\in H^m_G(M) \, .
$$

Let us now write $H^*_G(*)=\Ff_2[t]$, where $t=e(L)$,
so that $H^*_G(M)$ is the polynomial ring $H^*(M)[t]$. 
The identity says algebraically that
$$
a(t)\cdot t(t^d+w_1(\tau M)t^{d-1}+\ldots +w_d(\tau M))
=t^m\in H^*(M)[t],
$$
where $a(t)$ is a polynomial of degree $m-d-1$.
Multiplying by $\sum_{i\geq 0} t^{-i}w_i(-\tau M)$
in $H^*(M)[t,t^{-1}]$, we see that
$a(t)=t^{m-d-1}\sum_{i\geq 0} t^{-i}w_i(-\tau M)$.
Hence $w_i(-\tau M)=0$ for $i\geq m-d$. 
\end{proof}
As a corollary we obtain the classical criterion on normal
Stiefel-Whitney classes for the existence of an embedding
of $M$ in $\Rr^m$.
The usual proof for a smooth embedding involves showing
that the $\Ff_2$-Euler class of the normal bundle is zero
and requires less effort;
see Appendix \ref{atiyah} for a similar argument.

\Def{\label{D}
From now on we shall write $D$ for the greatest integer
$n\geq d$ such that $w_{n-d}(-\tau M)$ is non-zero.
Clearly $w_i(-\tau M)=0$ for $i>d$.
In fact, $e(\lambda )^{2d}=0$, because 
$H^{2d}(Y)$, Poincar\'e dual to $H_0(Y,\partial Y)$
for the $2d$-manifold $Y$,
is zero if $M$ is connected.
Thus we have the bounds: $d\leq D < 2d$.
}
\Prop{\label{app}
Suppose that $f: M\to \Rr^m$ is a continuous map.
If $m\leq D$, 
then there exist distinct points $x\not=y\in M$
such that $f(x)=f(y)$.
}
\begin{proof}
In the setting of Lemma \ref{embed},
the $G$-map $X\to L\otimes\Rr^m$:
$(x,y)\mapsto f(x)-f(y)$
gives, on orbit spaces, a section $s$ of $\lambda\otimes\Rr^m$.
Since $e(\lambda \otimes\Rr^m)=e(\lambda )^m$ is
non-zero by Lemma \ref{embed}, 
the section $s$ must have a zero $[x,y]$,
where $(x,y)\in X$.
\end{proof}
\Ex{Let $M$ be the real projective space of dimension
$d$. Let $r\geq 0$ be the integer such that $2^r \leq d<2^{r+1}$.
Then $D=2^{r+1}-1$.
}
\begin{proof}
Let $H$ be the Hopf line bundle over the projective space.
Then $\Rr\oplus \tau M\iso (d+1)H$.
In $H^*(M)=\Ff_2[T]/(T^{d+1})$,
where $w_1(H)$ is the class of $T$, the Stiefel-Whitney
classes of $-\tau M$ are determined by
$$
\sum_{i=1}^d w_i(-\tau M)T^i \equiv (1+T)^{-d-1}
\equiv (1+T)^{2^{r+1}-d-1}\quad ({\rm mod}\, T^{d+1}),
$$
because $2^{r+1}\geq d+1$. Since $2^{r+1}-d-1\leq d$,
the result is clear.
\end{proof}
Peterson's result \cite{fp} that a real projective space
of dimension $d=2^r$ does not embed in $\Rr^{2d-1}$ thus follows from
Proposition \ref{app}.
\Lem{\label{parameter}
With the hypotheses of Lemma \ref{embed},
suppose that $\xi$ is a 
real vector bundle of dimension $k$ over a compact ENR
$B$ and that $s$ is a section of the pullback of $\xi$ to
$B\times M$. If $w_{k-m}(\xi )\not=0$ for some $m\leq D$,
there is a point $b\in B$ and distinct points
$x\not= y\in M$ such that
$s(b,x)=s(b,y)\in\xi_b$.
}
\begin{proof}
The section $s$ determines a $G$-equivariant section 
of $L\otimes \xi$ over
$B\times X$ taking the value $s(b,x)-s(b,y)$ at $(b,(x,y))$,
and an associated section $\sigma$ of 
$\lambda\otimes \xi$ over $B\times Y$.
We have 
$$
e(\lambda \otimes\xi )=\sum_{i=0}^k w_{k-i}(\xi )e(\lambda )^i
\in H^k(B\times Y)=\bigoplus_{i=0}^kH^{k-i}(B)\otimes H^i(Y). 
$$
So, since $e(\lambda)^m\not=0$, by
Lemma \ref{embed}, $e(\lambda \otimes\xi )\not=0$
and the section $\sigma$ must have a zero $(b,[x,y])$.
Then $s(b,x)=s(b,y)$.
\end{proof}
\Sect{Simplicial complexes}
Suppose that we are given a triangulation of the closed
$d$-manifold $M$:
$V$ is a finite set, $S$ is a set of non-empty subsets
of $V$ such that, if $\emptyset\not=J\subseteq I$
and $I\in S$, then $J\in S$, and we have a
homeomorphism from $M$ to the space $K$ of all continuous maps
$x : V\to [0,1]$ such that $\sum_{v\in V}x(v)=1$ and
the support of $x$,
$\supp (x) =\{ v\in V\st x(v)>0\}$, is in $S$.
\Lem{\label{KM}
For given $m\geq 1$, let $X$ be the subspace
$$
\{ (x,y)\in K\times K \st \supp (x)\cap\supp (y)=\emptyset,
\, \# \supp (x)-1+\#\supp (y)-1 \leq m\}
$$
of $K\times K$ with the free action of $G=\Zz /2$ interchanging
the factors.

If $m\leq D$, 
then $e(\lambda )^m\in H^m(Y)$ is non-zero.
}
\begin{proof}
Following Kishimoto and Matsushita \cite[Lemma 3.3]{kish},
we write down a $G$-equivariant map
$$
\{ (x,y)\in K\times K \st x\not=y\}
\to \tilde X=
\{ (x,y)\in K\times K\st\supp (x)\cap\supp (y)=\emptyset\}:
$$
$$
(x,y)\mapsto (\alpha (x,y),\alpha (y,x)),
$$
where
$$
\alpha (x,y)(v) = \max\{ x(v)-y(v),0\}/
(\sum_{u\in V} \max\{ x(u)-y(u),0\} )\, .
$$
Notice that, if $x\not=y$, then there must be some 
$u\in V$ such that $x(u)>y(u)$. Now
$\supp (\alpha (x,y))=\{ v\in V\st x(v)>y(v)\}$.
So $\supp (\alpha (x,y))\cap\supp (\alpha (y,x))=\emptyset$.

We deduce from Lemma \ref{embed}
that the Euler class of $\tilde\lambda =\tilde X\times_GL$
over $\tilde Y=\tilde X/G$ satisfies
$e(\tilde\lambda )^m\not=0$ and hence that
its restriction $e(\lambda )^m$ to the
$m$-skeleton $Y$ (that is, $X/G$ in our standard notation) of 
$$
\tilde Y=\{ (x,y)\in K\times K\st \supp (x)\cap\supp (y)=\emptyset\} /G
$$
is non-zero.
\end{proof}
\Rem{\label{sphere}
Suppose that $\# V=d+2$ and $S$ is the set of all non-empty subsets
$I$ with cardinality $\# I\leq d+1$. Thus $K$ is the boundary of the simplex
with vertex set $V$ and is homeomorphic to the unit sphere $S(E)=M$
in the vector space $E$ of maps $z: V\to\Rr$
such that $\sum_{v\in V}z(v)=0$. 

We have an evident $G$-homeomorphism
$$
\tilde X \to S(L\otimes E)\quad :\quad
(x,y)\mapsto (x-y)/ \| x-y\|\, .
$$
Thus $\tilde\lambda$ corresponds to the Hopf line bundle
over the real projective space $P(E)$ of dimension $d$, and hence
we know, simply from the computation of the cohomology of the
projective space, that
$e(\hat\lambda )^d\in H^d(\hat Y)$
is non-zero,
without appealing to Lemma \ref{embed} for $m=d$ and $X=\tilde X$.
}
\Prop{\label{kish}
{\rm (\cite[Theorem 3.4]{kish}).}
Suppose that $f: K\to \Rr^m$ is a continuous map.
If $m\leq D$, 
then there exist points $x,\, y\in K$
such that $f(x)=f(y)$ with $\supp (x)\cap\supp (y)=\emptyset$
and $\#\supp (x)-1+\#\supp (y)-1\leq m$.
}
For the boundary of a simplex as in Remark \ref{sphere} this
reduces to the topological Radon theorem \cite{radon}.
\begin{proof}
The assertion follows from Lemma \ref{KM} just as
Proposition \ref{app} follows from Lemma \ref{embed}.
\end{proof}

\Lem{\label{flores}
For positive integers $m$ and $r$,
suppose that
$R_1,\ldots ,R_r$ are subsets of $S$
with the property that,
for each $j=1, \ldots ,r$,
if $I$ and $J$ are disjoint
subsets of $V$ in $S$ with $\# I-1 +\# J-1\leq m+r$, 
then either $I$ or $J$ is in $R_j$.

Write $R =\bigcap_{j=1}^r R_j$ and $A\subseteq K$ for the subcomplex
of maps $x\in K$ such that $\supp (x)\in R$.
Let $X$ be the subspace
$$
\{ (x,y)\in A\times A \st \supp (x)\cap\supp (y)=\emptyset\}
$$
of $A\times A$ with the free action of $G=\Zz /2$ interchanging
the factors.

If $m+r\leq D$, 
then $e(\lambda )^m\in H^m(Y)$ is non-zero.
}
\begin{proof}
We follow the now standard van Kampen-Flores argument
(as used, implicitly, in the proof of \cite[Lemma 4.2]{frick3}).

Define, for $j=1,\ldots ,r$,
a continuous function $\mu_j : K\to\Rr$ by
$$
\mu_j(x)=\min\{ \sum_{v\in V-I} x(v)\st I\in R_j\sqcup \{\emptyset\}\}.
$$
Then $\mu_j(x)=0$ if and only if $\supp (x)\in R_j$.

Write $\hat X$ for the $G$-subspace
$$
\{ (x,y)\in K\times K \st \supp (x)\cap\supp (y)=\emptyset,\,
\# \supp (x)+\# \supp (y)\leq m+r+2\}
$$
of $K\times K$.
The $G$-map $\hat X\to L\otimes\Rr^r$:
$(x,y)\mapsto (\mu_j(x)-\mu_j(y))_{j=1}^r$
determines a section $s$ of $\hat\lambda\otimes \Rr^r$
(where, in our usual notation, $\tilde\lambda =\tilde X\times_G L$).
If $\mu_j (x)=\mu_j (y)$, $\supp (x)\cap\supp (y)=\emptyset$
and $\#\supp (x)+\#\supp (y)\leq m+r+2$, 
then one of $\supp (x)$ and $\supp (y)$ is in $R_j$.
Hence one of $\mu_j (x)$, $\mu_j (y)$ is zero,
and so both are. Thus the zero-set of the section $s$ is precisely
$Y=X/G$.

Since $m+r\leq D$, we have $e(\hat \lambda )^{m+r}\not=0$.
It follows, by a standard argument recalled in
Appendix \ref{zero}, that $e(\lambda )^m\not=0$. 
\end{proof}
\Rem{Assuming that $R\not=S$, let $\Gamma$ be the set of minimal
elements, with respect to $\subseteq$, in the complement $S-R$.
The subsets $C_j=\{ I'\in\Gamma \st I'\notin R_j\}$ have the 
properties:

(i)
$\Gamma =\bigcup_{i=1}^r C_j$;

(ii)
if $I',\, J'\in C_j$, then either $\# I' +\# J'>m+r+2$ or 
$I'\cap J'\not=\emptyset$.
\par\noindent
Conversely, if subsets $C_j\subseteq\Gamma$ satisfy (i) and (ii),
then the subsets $R_j=\{ I\in S\st I'\not\subseteq I$
for all $I'\in C_j\}$ of $S$ satisfy the conditions 
required in Lemma \ref{flores}.
}
\Prop{\label{kish2}
With the hypotheses of Lemma \ref{flores},
suppose that $f: A\to\Rr^m$ is a continuous map.
If $m+r\leq D$, 
then there exist $x,\, y\in A$ with
disjoint support such that $f(x)=f(y)$.
}
\begin{proof}
This is a now routine deduction from Lemma \ref{flores}.
\end{proof}
For $r=1$, the condition is satisfied by
$R_1=\{ I\in S \st 2\# I \leq m+3\}$, 
and, if $m=2q$ is even, $A$ is the $q$-skeleton $K^{(q)}$.
Thus we obtain the theorem
\cite[Theorem 1.2]{kish} of Kishimoto and Matsushita.
\Cor{Suppose that $2q+1\leq D$.
Then, for any map $f : K^{(q)}\to \Rr^{2q}$ on the $q$-skeleton
of $K$, there exist $x,\, y\in K^{(q)}$ with disjoint support
such that $f(x)=f(y)$.
\qed
}
The condition is clearly satisfied if $d=2q+1$,
or if $d=2q$ and some Stiefel-Whitney class $w_i(\tau M)$
with $i\geq 1$ is non-zero.
For the boundary of a simplex as in Remark \ref{sphere}
with $A$ a skeleton of the simplex,
Proposition \ref{kish2} is the classical van Kampen-Flores theorem
\cite{flores, kampen}.
\Rem{Continuing Remark \ref{sphere},
suppose that $V=\bigsqcup_{j=1}^{r} V_j$ is a disjoint union
of $r$ non-empty subsets.
Take $R_j=\{ I\in S \st 2\# (I\cap V_j) < \# V_j\}$;
Clearly, if $I$ and $J$ are disjoint,
then either $I$ or $J$ is in $R_j$.
The simplicial complex $A$ is the $j$-fold join of skeleta
of $\Delta (V_j)$.
We can apply Proposition \ref{kish2} if
$m+r\leq D=d$.   
For example, if $\# V_j=3$, $m=2q$, $r=q+1$,
as considered in \cite[Corollary 1.3]{frick3}.
}

\Lem{\label{app2}
With the hypotheses of Lemma \ref{flores}, 
suppose that $\xi$ is a $k$-dimensional
real vector bundle over a compact ENR
$B$ and that $s$ is a section of the pullback of $\xi$ to
$B\times A$. If 
$w_{k-m}(\xi )\not=0$ for some $m\leq D-r$,
there is a point $b\in B$ and
$x,\, y\in A$ with disjoint support such that
$s(b,x)=s(b,y)\in\xi_b$.
}
\begin{proof}
This follows like the deduction of Lemma \ref{parameter}
from Lemma \ref{embed}.
\end{proof}
As an application we have the following result of Frick and Harrison.
\Prop{Let $T : \Rr^{m+1}\to\Rr^{m+1}$ be a reflection in a hyperplane.
With the hypotheses of Lemma \ref{flores}, 
suppose that $m+r\leq D$.
Then, for any homotopy
$h_t : A \to \Rr^{m+1}$,
$0\leq t\leq 1$ 
from a map $f=h_0$ to its `mirror image'
$T\comp f=h_1$, 
there exist some $t\in [0,1]$ and $x,\, y\in A$
with disjoint support such that $h_t(x)=h_t(y)$.
}
\begin{proof}
Suppose that $T$ is the reflection in the hyperplane $\Rr^m\oplus 0$.
We take $B$ to be the projective line $S(L\oplus L)/G$
and $\xi$ to be direct sum 
$S(L\oplus L)\times_G (\Rr^m\oplus L)=\Rr^m\oplus H$,
where $H$ is the Hopf line bundle.
So $k=m+1$, and $w_{k-m}(\xi )=w_1(H)\not=0$.

The homotopy defines a $G$-map
$F: S(L\oplus L)\times A \to \Rr^m\oplus L$ 
by 
$$
F((\cos (\pi t),\sin (\pi t)),x)= h_t(x), \quad0\leq t\leq 1,
$$
and this determines a section of $\xi =\Rr^m\oplus H$.
The conclusion of Lemma \ref{app2}
asserts that there exist $t\in [0,1]$ and $x,\, y\in A$
with disjoint support such that $h_t(x)=h_t(y)$.
\end{proof}
\Cor{Let $T : \Rr^{2q+1}\to\Rr^{2q+1}$ be a reflection in a hyperplane.

Suppose that 
$2q+1\leq D$.
Then, for any homotopy
$h_t : K^{(q)} \to \Rr^{2q+1}$,
$0\leq t\leq 1$ 
from a map $f=h_0$ to its `mirror image'
$T\comp f=h_1$, 
there exist some $t\in [0,1]$ and $x,\, y\in K^{(q)}$
with disjoint support such that $h_t(x)=h_t(y)$.
\qed
}

\Ex{{\rm (\cite[Theorem 1.5]{frick3}).}
More generally,
under the hypotheses of Lemma \ref{flores}, 
suppose given integers 
$0\leq l\leq m\leq k$, such that $m+r\leq D$ and
$\binom{k-l}{k-m}$ is odd.
Then for every $G$-map 
$F: S((k-m+1)L)\times A\to\Rr^l\oplus L^{k-l}$
there exist $e\in S((k-m+1)L)$ and $x,\, y\in A$
with disjoint support
such that $F(e,x)=F(e,y)$.
}
\begin{proof}
Take $B=S((k-m+1)L)/G$,
a real projective space of dimension $k-m$, and 
$\xi=S((k-m+1)L)\times_G(\Rr^l\oplus L^{k-l})= \Rr^l\oplus (k-l) H$,
the direct sum of a trivial bundle
and $k-l$ copies of the Hopf line bundle $H$.
We have $w_{k-m}(\xi )\not=0$, since $\binom{k-l}{k-m}$ is odd.
\end{proof}

\Sect{Discussion}
\Lem{In the framework of Lemma \ref{embed},
suppose that $m>d$ is an integer such that
$w_i(-\tau M)=0$ for all $i\geq m-d$.
Then $e(\lambda )^m=0$.
}
\begin{proof}[Outline proof]
Reversing the argument in the proof of Lemma \ref{embed},
we obtain an element $a\in H^{m-d-1}_G(M)$ such that
$\Delta_*(t\cdot\Delta_!(a))=t^m\in H^m_G(M)$.
The cohomology ring $H^*_G(M\times M)$ can be described
using the Steenrod squaring operation
$P^2: H^i(M)\to H^{2i}_G(M\times M)$ and the
induction homomorphism $\ind : H^i(M\times M) \to H^i_G(M\times M)$.
As an $\Ff_2[t]$-module, $H^*_G(M\times M)$ is generated
by the classes $P^2(u)$ and
$\ind (v\otimes w)$, for $u,\,v,\, w\in H^*(M)$.
The image of $\ind$ is annihilated by multiplication by
$t$;  the composition $t\cdot P^2 :H^i(M)\to H^{2i+1}_G(M\times M)$
is linear; and the restriction of $P^2(u)$ to the
diagonal defines the Steenrod squares: 
$$
\Delta^*(P^2(u))=\sum_{j=0}^i t^{i-j}\Sq^{j}u
\in H^*(M)[t]
\text{\quad for $u\in H^i(M)$.}
$$
Thus we can write $t\cdot\Delta_!(a)-t^m
=\sum_{2i\leq m} t^{m-2i}P^2(u_i)$, with $u_i\in H^i(M)$.
So $\sum_{2i\leq m} t^{m-2i}\Delta^*(P^2(u_i))=0\in H^*(M)[t]$,
that is, $\sum_{2i\leq m, \, j\leq i} t^{m-i-j}\Sq^ju_i=0$.
Hence $u_i=0$ for all $i$. We conclude that 
$t^m =\Delta_!(t\cdot a)\in H^m_G(M\times M)$
and so $e(\lambda )^m=0\in H^m(Y)$.
\end{proof}

Thus the integer $D$ in Definition \ref{D} is the greatest integer
$n$ such that $e(\lambda )^n\not=0$.
We have used only the Euler class in $\Ff_2$-cohomology.
In principle one can use any multiplicative cohomology theory.
For example, in \cite{mark} Mahowald showed, in effect, that
for $M$ a real projective space of dimension $d=2^r+1$,
$r\geq 1$, the
integral cohomology class $e_\Zz (\lambda )$ (with twisted
coefficients) satisfies $e_\Zz (\lambda )^{2d-2}\not=0$,
whereas $D=2d-3$,
with the consequence that $M$ does not embed in $\Rr^{2d-2}$.

Writing $\omega^*$ for stable cohomotopy theory, one
has a universal stable cohomotopy Euler class\footnote{See, for example, \cite[II.4]{FHT}.}
that we shall write as
$$
\gamma (\lambda )\in \omega^0_G(Y;\, -\lambda )
$$
in the stable cohomotopy of the Thom space of 
the  virtual bundle $-\lambda$ over $Y$. 
If one replaces $D$ by the greatest integer $n\geq d$ such that
$\gamma (\lambda )^n$ is non-zero, the 
Propositions \ref{app}, \ref{kish} and \ref{kish2}
(but not Lemmas \ref{parameter} and \ref{app2})
are still valid, with unchanged proofs.
But, of course, the determination of this integer 
for stable cohomotopy is likely to be very difficult. 

The $2$-primary torsion 
class $\gamma (\lambda )^n\in\omega^0(Y;\, -n\lambda )$ is an obstruction to the existence of a continuous injective map,
that is, a topological embedding, $M\into \Rr^n$.
Haefliger \cite{haefliger} showed that, in a stable range
$3d<2(n-1)$, this is the precise obstruction to the existence of a 
smooth embedding $M\into\Rr^n$.

For smaller $n$, there are other obstructions 
(not only at the prime $2$) to the existence
of a smooth, rather than topological, embedding.
Some remarks on such obstructions in integral cohomology
and topological periodic $K$-theory can be found in
Appendix \ref{atiyah}.
\begin{appendix}
\Sect{\label{atiyah} Smooth embeddings}
Suppose that
we have a smooth embedding $f: M\into \Rr^n$ with normal bundle
$\nu$. 
We shall formulate the classical restriction on Pontryagin classes
in terms of complexification,
giving complex vector bundles the standard orientation as real
bundles.
\Prop{If there is a smooth embedding $M\into\Rr^n$, then
$c_i(-\Cc\otimes\tau M)=0\in H^{2i}(M;\,\Zz )$ for
$i\geq n-d$.
}
\begin{proof}
The Chern classes $c_i(\Cc\otimes\nu )=c_i(-\Cc\otimes\tau M)$
of the complexified normal bundle $\Cc\otimes\nu$ of complex dimension
$n-d$ are zero if $i>n-d$.
This is a condition for existence of an immersion.
The condition that $f$ be an embedding implies, further, 
that the $\Zz$-cohomology Euler class
$e_\Zz (\Cc\otimes\nu )=0$. This can be seen as follows.

Consider a tubular neighbourhood $D(\nu )\into B_R(\Rr^n)$
into the open ball of sufficiently large radius $R$ in $\Rr^n$.
We have a commutative diagram
$$
\begin{matrix}
H^{2(n-d)}(D_R(\Rr^n),D_R(\Rr^n)-B(\nu );\,\Zz )&\Rarr{\iso}{}&
H^{2(n-d)}(D(\nu ),S(\nu );\,\Zz ) \\
\Darr{}{}&&\Darr{}{}\\
H^{2(n-d)}(D_R(\Rr^n);\,\Zz )=H^{2(n-d)}(*;\,\Zz )
&\Rarr{}{} & H^{2(n-d)}(D(\nu );\,\Zz )
=H^{2(n-d)}(M;\,\Zz )\\
\Darr{}{}&&\Darr{}{}\\
H^{2(n-d)}(D_R(\Rr^n)-B(\nu );\,\Zz )
&\Rarr{}{}& H^{2(n-d)}(S(\nu );\,\Zz )
\end{matrix}
$$
of exact sequences of pairs.
The lift of $e_\Zz (\Cc\otimes \nu )$ from $M$ to the sphere bundle
$S(\nu )$ is zero, because
the pullback of $\nu$ to $S(\nu )$ has a tautological
nowhere-zero section.
Thus 
$e_\Zz (\Cc\otimes\nu )$ lies in the image of $H^{2(n-d)}(*;\,\Zz )$
under the homomorphism induced by the projection $M\to *$.
But $e_\Zz (\Cc\otimes \nu )$ restricts to $0$ at any point $*\in M$,
since $n-d>0$. Hence $e_\Zz (\Cc\otimes\nu )=0$.

(The argument has been given without exploiting the fact that
$H^{2(n-d)}(*;\,\Zz )=0$ and so carries through with
only notational changes in periodic $K$-theory, too.)
\end{proof}
\Rem{Notice that the condition $c_i(-\Cc\otimes\tau M)=0$
for $i\geq n-d$ is automatically satisfied if
$2(n-d)>d$, that is, $3d<2n$.
}
We can rewrite the condition in a form closer to that
used in Lemma \ref{embed}.
Let us write $E$ for the representation $\Cc$ of the
circle group $\Tt$ of complex numbers of modulus $1$ acting
by multiplication. 

The complexification of the derivative $(df)_x : \tau_xM\into\Rr^n$ 
at $x\in M$, is an injective
$\Tt$-equivariant linear map $E\otimes \tau_xM \to E\otimes \Rr^n$.
The Euler class $e_\Zz (E)^n\in H^{2n}_\Tt (M;\,\Zz )$ 
in Borel cohomology with $\Zz$-coefficients factors as
$$
a\cdot e_\Zz (E)\cdot e_\Zz (E\otimes\tau M) =e_\Zz (E)^n\, .
$$
Writing $t=e_\Zz (E)$, we have
$$
a(t)\cdot t
\cdot (t^d+c_1(\Cc\otimes\tau M)t^{d-1}+
\ldots +c_d(\Cc\otimes\tau M))=t^n,
$$
in $H^*_\Tt (M;\,\Zz )=H^*(M;\,\Zz )[t]$,
where $a(t)$ is a polynomial of degree $n-d-1$.

In \cite{mfa} Atiyah used the same argument for periodic
complex $K$-theory.
The $K$-theory Euler class of a complex vector bundle
$\xi$ of dimension $k$ over a compact ENR $B$,
as defined by the standard Bott class in $K^{2k}(D(\xi ),S(\xi ))$,
can be written as $v^{-k}\sum_{i=0}^k(-1)^i
\lambda^i\xi$ in $K^{2k}(B)$, 
where $v\in K^{-2}(*)=\Zz v$ is the Bott generator
and $\lambda^i\xi$ is the class of the $i$th exterior power of $\xi$.
$K$-theory Chern classes $c_i^K(\xi )\in K^{2i}(B)$ can be defined
using the $\Tt$-equivariant $K$-theory $K^*_\Tt (B)=K^0(B)[z,z^{-1}]$,
where $z=[E]\in K^0_\Tt (*)$, by
$$
c_0^K(\xi )(1-z)^k+\ldots + v^ic_i^K(\xi )(1-z)^{k-i}+\ldots +
v^kc_k^K(\xi )=
\sum_{i=0}^k(-z)^i\lambda^i\xi\, .
$$
Thus $c_0^K(\xi )=\lambda^k\xi$, which is a unit in $K^0(B)$,
$c_k^K(\xi )=e_K(\xi )$ is the $K$-Euler class,
and $c_0^K(\xi)+vc_1^H(\xi )+\ldots +v^kc_k^K(\xi )=1$.
If $\xi$ is trivial, $c_0^K(\xi )=1$ and $c_i^K(\xi )=0$ for $i>1$.
And the Chern classes are multiplicative: if $\eta$ is another
complex bundle over $B$, then
$c_i^K(\xi\oplus\eta )=\sum_{r+s=i} c_r^K(\xi )c_s^K(\eta )$.

It follows, just as in cohomology,
that $c_i^K(-\tau M)=0$ for $i\geq n-d$
(and again this is automatically true if $3d<2n$, because
the Chern classes lift to connective $K$-theory).

In the $\Tt$-equivariant group
$K^0_\Tt (M)=K^0(M)[z,z^{-1}]$ we have
$$
a(z)\cdot (1-z)
\cdot \sum_{i=0}^d (-z)^i\lambda^i(\Cc\otimes\tau M )=
(1-z)^n,
$$
where $a(z)$ is a polynomial of degree $n-d-1$.
\Rem{In \cite{mfa} the condition is stated in terms
of the Grothendieck operations $\gamma^i$, which are related to
the Chern classes by the formal power series identity
$$
\sum_{i\geq 0}T^i\gamma^i(\xi -k) =\sum_{i=0}^k (1-T)^iv^ic^K_i(\xi )
\in K^0(B)[[T]],
$$
as $\gamma^i(\Cc^k-\Cc\otimes\tau M)=0$ for $i\geq n-d$.
By using Real $K\Rr$-theory one obtains Atiyah's condition
\cite[Theorem 4.3]{mfa}:
$\gamma^i(\Rr^d-\tau M)=0\in \KO^0(M)$ for $i\geq n-d$.
}
\Sect{\label{zero} A lemma on the zero-set of a section}
\Lem{
Let $\xi$ and $\eta$ be real vector bundles of dimension $m$ and $n$
over
a compact ENR $B$, and suppose that a closed sub-ENR
$Z\subseteq B$ is the zero-set of a section of
$\eta$. 

\par\noindent {\rm (a)}
If the $\Ff_2$-cohomology Euler class
$e(\xi\oplus\eta )\in H^{m+n}(B)$ is non-zero,
then the Euler class $e(\xi\, |\, Z)\in H^m(Z)$
of the restriction of $\xi$ to $Z$ is non-zero.

\par\noindent {\rm (b)}
If the stable cohomotopy Euler class
$\gamma (\xi\oplus\eta )\in \omega^0(B;\, -(\xi\oplus\eta ))$ 
is non-zero,
then $\gamma (\xi\, |\, Z)\in \omega^0(Z;\, -\xi\, |\, Z)$
is non-zero.
}
\begin{proof}
Let $z$ be a section of $\eta$ with zero-set $Z$.

\par\noindent
(a).
A rather more general result is proved in 
\cite[Proposition 2.7]{borsuk}. We give a proof here
using singular cohomology.

Let a class $a\in H^m(V)$ restrict to zero in $H^m(Z)$.
Since $Z$ is assumed to be a closed sub-ENR, there is an open
neighbourhood $U$ of $Z$ in $B$ and a deformation
retraction $r_t: U\to B$, $0\leq t\leq 1$, such that
$r_t(x)=x$ for $x\in Z$, $r_0(x)=x$ for all $x\in U$,
$r_1(U)=Z$. So $a$ restricts to zero in $H^m(U)$.
The Euler class $e(\eta )$ restricts to zero in 
$H^n(B-Z)$.
Since $B$ is the union of the two open sets
$U$ and $V=B-Z$ and the cup product lifts to
$H^m(B,U)\otimes H^n(B,V)\to H^{m+n}(B,U\cup V)$,
the product $a\cdot e(\eta )$ must be zero. 

As $e(\xi\oplus \eta )=e(\xi )\cdot e(\eta )$,
we conclude that $e(\xi )\not=0$.

Alternatively, 
if $B$ is a simplicial complex and $Z$ is a subcomplex, we can find
subcomplexes $P$ and $Q$ of a subdivision of $B$ such that
$P\subseteq U$, $Q\subseteq V$, $P\cup Q=B$, and use the product
$H^m(B,P)\otimes H^n(B,Q)\to H^{m+n}(B,P\cup Q)$.

\par\noindent
(b). 
Suppose that $\gamma (\xi\, |\, Z)=0$. Then, for sufficiently large
$k\geq 1$,
there is a nowhere-zero
section $s$ of the pullback of $\Rr^k\oplus\xi$ to
$D(\Rr^k)\times Z$ such that, for $(v,x)\in S(\Rr^k)\times Z$
$s(v,x)=(v,0)$. This section $s$ extends to a section $\tilde s$
over $D(\Rr^k)\times B$ such that $\tilde s(v,x)=(v,0)$ for
$(v,x)\in S(\Rr^k)\times B$.
The section $(v,x)\mapsto (\tilde s(v,x), (1-\| v\|)z(x))$
of $\Rr^k\oplus \xi\oplus \eta$
over $D(\Rr^k)\times B$ is nowhere zero and takes the
value $(v,0,0)$ for $v\in S(\Rr^k)$.
So $\gamma (\xi\oplus\eta )$ must be zero.
\end{proof}
\Rem{The simpler result that 
\par\noindent
{\it if either 
{\rm (a)} $e(\xi\oplus\eta )\not=0$ or 
{\rm (b)} $\gamma (\xi\oplus\eta )\not=0$,
then any section of $\xi\, |\, Z$ has a zero} 
\par\noindent
is easily proved.

For any section $s$ of $\xi\, |\, Z$ extends to a section
$\tilde s$ of $\xi$. If $z$ is a section of $\eta$ with zero-set $Z$,
then the zero-set of the section $(\tilde s,z)$ of $\xi\oplus\eta$
is precisely the zero-set of the given section $s$.
(Of course, (a) implies (b);
because $e(\xi\oplus\eta )$ is the Hurewicz image of 
$\gamma (\xi\oplus\eta )$.)
}
\end{appendix}

\end{document}